\documentclass[titlepage,12pt]{article}

\usepackage[utf8]{inputenc}
\usepackage[T1]{fontenc}
\usepackage[english]{babel}
\usepackage[a4paper,top=2cm,bottom=2cm,left=3cm,right=3cm,marginparwidth=1.75cm]{geometry}
\usepackage{graphicx}
\graphicspath{{figures/}}
\usepackage{subcaption}
\usepackage{caption}
\usepackage[colorlinks=true,allcolors=blue,hypertexnames=false]{hyperref}
\usepackage{listings}
\usepackage{spverbatim}
\usepackage[section]{placeins}
\usepackage{amsmath}
\usepackage{amsthm}
\usepackage{mathtools}
\usepackage{amsfonts}
\usepackage{longtable}
\usepackage{array}

\author{Jakub Pączka \\
Agnieszka Przeździecka \\
Piotr Tyrakowski}

\title{%
    Forbidden Substructures Preserved by the Blob Graph Operation
    \thanks{The original version of this paper was prepared as part of the course ``Research Project: Structural Graph Theory'' at the Faculty of Mathematics and Information Science, Warsaw University of Technology, under the supervision of Professor Paweł Rzążewski during the 2024/2025 academic year.}
}

\date{October 2026}

\newtheorem{theorem}{Theorem}[section]
\newtheorem{definition}{Definition}
\newtheorem{lemma}{Lemma}[theorem]
\newtheorem{corollary}{Corollary}[theorem]

\begin{document}
\maketitle

\begin{abstract}
    The blob graph \(G^\circ\) of a graph \(G\) has a vertex for each non-empty connected vertex subset of \(G\), with two distinct subsets adjacent whenever they intersect or are joined by an edge of \(G\). This construction has been used in reductions from induced subgraph optimisation problems to maximum-weight independent set, where the preservation of forbidden substructures is essential. We completely characterise the graphs \(H\) for which excluding \(H\) from \(G\) guarantees that \(G^\circ\) also excludes \(H\), under six containment relations. For subgraphs, minors, and topological minors, preservation holds exactly when \(H\) has at most two vertices. For induced subgraphs, it holds exactly when \(H\) is a linear forest. For induced topological minors, it holds exactly when \(H\) is triangle-free and has maximum degree at most two. For induced minors, it holds exactly when \(H\) is triangle-free and every neighbour of a vertex of degree at least three has degree exactly two. These classifications unify and extend preservation results previously established for particular forbidden graphs in algorithmic applications.
\end{abstract}

\section{Introduction}

The goal of this paper is to study the properties of blob graphs. The blob graphs are graphs created from a graph \(G\), whose vertices are connected subgraphs of \(G\). Two vertices of the blob graph are adjacent if the corresponding subgraphs of \(G\) intersect, or if there is an edge of \(G\) with one endpoint in one subgraph and the other endpoint in the second subgraph. They are properly defined in Definition~\ref{def:blob}. The blob graph of a graph \(G\) is denoted by \(G^\circ\).

\paragraph{Motivation.} The blob graph definition was introduced as a tool for reducing the problem of finding a maximum-weight induced subgraph of \(G\) whose components lie in a fixed hereditary graph class to the Max Weight Independent Set problem on an induced subgraph of \(G^\circ\)\cite{gartland2021findinglargeinducedsparse}. This reduction and its later uses~\cite{paesani2022feedbackvertexsetcycle}\cite{bonnet2026qptasmwisfindinglarge} rely on the fact that the certain forbidden-substructure properties of \(G\), transfer to \(G^\circ\): excluding a long induced cycle~\cite{gartland2021findinglargeinducedsparse}, excluding a fixed linear forest as an induced subgraph~\cite{paesani2022feedbackvertexsetcycle}, or excluding a disjoint union of long induced cycles as an induced minor~\cite{bonnet2026qptasmwisfindinglarge}. Each of these results was estabilished only for a specific graph \(H\) needed in the corresponding application. This suggests a natural question: for which graphs \(H\), and for which containment relation, does excluding \(H\) from \(G\) guarantee that \(G^\circ\) also excludes \(H\)? In this paper we answer this question completely.

\paragraph{Our results.} Our main contribution is a complete classification of the graphs \(H\) for which \(H\)-freeness is preserved under the blob graph operation, under each of the following containment relations: subgraphs, minors, topological minors, induced subgraphs, induced topological subgraphs and induced minors. For each of these notions we determine exactly for which graphs \(H\) the implication ``if \(G\) excludes \(H\), then \(G^\circ\) also excludes \(H\)'' holds for every graph \(G\).

Our characterization for induced subgraphs (Theorem~\ref{theo:podind}) is closely related to a fact used by Paesani, Paulusma and Rzążewski~\cite{paesani2022feedbackvertexsetcycle}, and our characterisation for induced minors (Theorem~\ref{theo:minind}), restricted to graphs \(H\) that are disjoint unions of long induced cycles, reintroduces results of Gartland et al.~\cite{gartland2021findinglargeinducedsparse} and of Bonnet et al.~\cite{bonnet2026qptasmwisfindinglarge}. The general statement of Theorem~\ref{theo:minind} extends these special cases to arbitrary triangle-free graph \(H\) in which every neighbour of a vertex of degree at least three has degree exactly two.

The paper is organized as follows. In section \ref{prelim} we introduce notation and preliminary definitions. Sections \ref{pod}--\ref{minind} treat subgraph, minor, topological minor, induced subgraph, induced topological minor and induced minor invariance, respectively.

\section{Preliminaries}\label{prelim}

For any graph \(G\), we denote its blob graph by \(G^\circ\) and define it as follows.

\begin{definition}\label{def:blob}
    \[
        G^\circ \coloneq (V^\circ,E^\circ),
    \]
    where
    \begin{align*}
        V^\circ & \coloneq \left\{ X\subseteq V(G): G[X]\text{ is connected} \right\}, \\
        E^\circ & \coloneq \left\{XY: X,Y\in V^\circ,\ X\neq Y,\
        \left(X\cap Y\neq\emptyset\ \lor\
        \exists x\in X,\exists y\in Y \text{ such that } xy\in E(G)\right)\right\}.
    \end{align*}
\end{definition}

In other words, the vertices of \(G^\circ\) are the non-empty connected vertex subsets of \(G\). Two such subsets are adjacent in \(G^\circ\) if they intersect, or if there is an edge of \(G\) with one endpoint in one subset and the other endpoint in the second subset.

\section{Subgraphs}\label{pod}

We first determine for which graphs \(H\) the following definition holds.

\begin{definition}[\textbf{S-invariant}]\label{def:pod}
    A graph \(H\) is called \textbf{S-invariant} if for every graph \(G\), the fact that \(G\) does not contain \(H\) as a subgraph implies that \(G^\circ\) does not contain \(H\) as a subgraph.
\end{definition}

\begin{theorem}\label{theo:pod}
    A graph \(H\) is S-invariant if and only if \(|V(H)|\leq 2\).
\end{theorem}

\begin{proof}[Proof of the forward implication]
    Suppose that a graph \(H\) is S-invariant, and assume that \(|V(H)|\geq 3\). Put \(n=|V(H)|\), and let \(G=P_{n-1}\) be the path with vertex set \(\{v_1,v_2,\dots,v_{n-1}\}\). Clearly, \(G\) does not contain \(H\) as a subgraph, because \(|V(G)|<|V(H)|\). We show that \(G^\circ\) contains \(H\) as a subgraph.

    Consider the set
    \[
        X\coloneq \bigl\{\{v_2\},\{v_1\},\{v_1,v_2\},\{v_1,v_2,v_3\},\dots,\{v_1,v_2,\dots,v_{n-1}\}\bigr\}.
    \]
    We claim that \(G^\circ[X]\) is a clique. The edge between \(\{v_1\}\) and \(\{v_2\}\) exists because \(v_1v_2\in E(G)\). For every other pair of vertices in \(X\), the corresponding subsets have non-empty intersection, containing at least one of \(v_1\) or \(v_2\). Hence all such pairs are adjacent in \(G^\circ\).

    Therefore \(G^\circ\) contains a clique on \(n\) vertices as a subgraph. Since \(K_n\) contains every graph on \(n\) vertices as a subgraph, it contains \(H\) in particular. Hence \(G^\circ\) contains \(H\) as a subgraph, a contradiction.

    Thus, no graph \(H\) with at least three vertices is S-invariant. This proves the forward implication.
\end{proof}

\begin{proof}[Proof of the reverse implication]
    Let \(H\) be a graph with at most two vertices. We consider all possibilities up to isomorphism.

    \begin{enumerate}
        \item \(H=(\{v\},\emptyset)\). Depending on the convention, either there is no graph \(G\) that does not contain a single vertex, or \(G\) is the empty graph. In the latter case, \(G^\circ\) is also empty, so it also does not contain \(H\) as a subgraph. In the remainder of the proof we adopt the convention that the empty graph does not exist.
        \item \(H=(\{v,u\},\emptyset)\). If \(G\) does not contain \(H\) as a subgraph, then \(G\) has exactly one vertex. Therefore \(G^\circ\) also has exactly one vertex, and hence does not contain \(H\) as a subgraph.
        \item \(H=(\{v,u\},\{uv\})\). If \(G\) does not contain \(H\) as a subgraph, then \(G\) contains no edge. Consequently, \(G^\circ\) also contains no edge: in an edgeless graph, the only connected non-empty vertex subsets are singletons, and distinct singletons are disjoint and non-adjacent. Thus \(G^\circ\) also does not contain \(H\) as a subgraph.
    \end{enumerate}

    Hence the reverse implication, and therefore the theorem, follows.
\end{proof}

\subsection{Generalisation to other non-induced substructures}

\begin{corollary}\label{con:pod}
    The following statements:
    \begin{itemize}
        \item for every graph \(G\), if \(G\) does not contain \(H\) as a minor, then \(G^\circ\) does not contain \(H\) as a minor;
        \item for every graph \(G\), if \(G\) does not contain \(H\) as a topological minor, then \(G^\circ\) does not contain \(H\) as a topological minor,
    \end{itemize}
    are true if and only if \(|V(H)|\leq 2\).
\end{corollary}

\begin{proof}
    The proof is analogous to the proof of Theorem~\ref{theo:pod}, so we do not repeat it in full. The key facts are the following:
    \begin{enumerate}
        \item No graph \(G\) contains a minor or a topological minor with more vertices than \(G\) itself.
        \item For every substructure that may be obtained using edge deletions, the complete graph \(K_n\) contains every possible such substructure on \(n\) vertices.
    \end{enumerate}
    These observations immediately imply Corollary~\ref{con:pod}.
\end{proof}

\section{Induced subgraphs}\label{podind}

\begin{definition}
    For graphs \(G\) and \(H\), we say that \(G\) is \(H\)-free if \(G\) does not contain \(H\) as an induced subgraph.
\end{definition}

We now determine for which graphs \(H\) the following definition holds.

\begin{definition}[\textbf{IS-invariant}]\label{def:podind}
    A graph \(H\) is called \textbf{IS-invariant} if for every graph \(G\), the fact that \(G\) is \(H\)-free implies that \(G^\circ\) is \(H\)-free.
\end{definition}

\begin{theorem}\label{theo:podind}
    A graph \(H\) is IS-invariant if and only if \(H\) is a subgraph of a path.
\end{theorem}

\begin{proof}[Proof of the forward implication]
    Suppose that a graph \(H\) is IS-invariant. Assume that \(H\) is not a subgraph of a path. Then \(H\) belongs to exactly one of the following families:
    \begin{itemize}
        \item graphs containing a cycle;
        \item forests that are not subgraphs of a path.
    \end{itemize}
    We treat the two cases separately. In each case, we construct a graph \(G\) such that \(G\) is \(H\)-free but \(G^\circ\) is not.

    \begin{lemma}\label{lem:podind1}
        For every graph \(H\) containing a cycle, there exists a graph \(G\) such that \(G\) is \(H\)-free, but \(G^\circ\) is not \(H\)-free.
    \end{lemma}

    \begin{proof}
        Assume that \(H\) contains a cycle. Let \(l\in\mathbb{N}\) be the length of a shortest cycle in \(H\). Clearly, \(l\geq 3\).

        Let \(G\) be the graph obtained from \(H\) by subdividing every edge into a path of length \(3\). Equivalently, every edge of \(H\) is replaced by a path in which two new vertices are inserted.

        We first show that \(G\) is \(H\)-free. Let \(C\) be any cycle in \(G\). Since \(G\) is a subdivision of \(H\), the cycle \(C\) corresponds to some cycle \(C'\) in \(H\), and \(|C|=3|C'|\). Hence \(|C|\geq 3l\). Thus \(G\) contains no cycle of length \(l\), and therefore cannot contain \(H\) as an induced subgraph.

        We now show that \(G^\circ\) is not \(H\)-free. Let \(X^\circ\) be the set of vertices of \(G^\circ\) corresponding to the maximal stars centred at the original vertices of \(H\). That is, for every \(u\in V(H)\), we take the subset of \(V(G)\) consisting of \(u\) and the subdivision vertices adjacent to \(u\). One can verify that an edge between two vertices of \(X^\circ\) exists if and only if their central vertices are adjacent in \(H\). Hence \(G^\circ[X^\circ]\simeq H\), so \(G^\circ\) contains \(H\) as an induced subgraph.

        Thus we have shown that no graph \(H\) containing a cycle is IS-invariant.
    \end{proof}

    \begin{figure}[h]
        \centering
        \begin{minipage}{.5\textwidth}
            \centering
            \includegraphics[width=.8\linewidth]{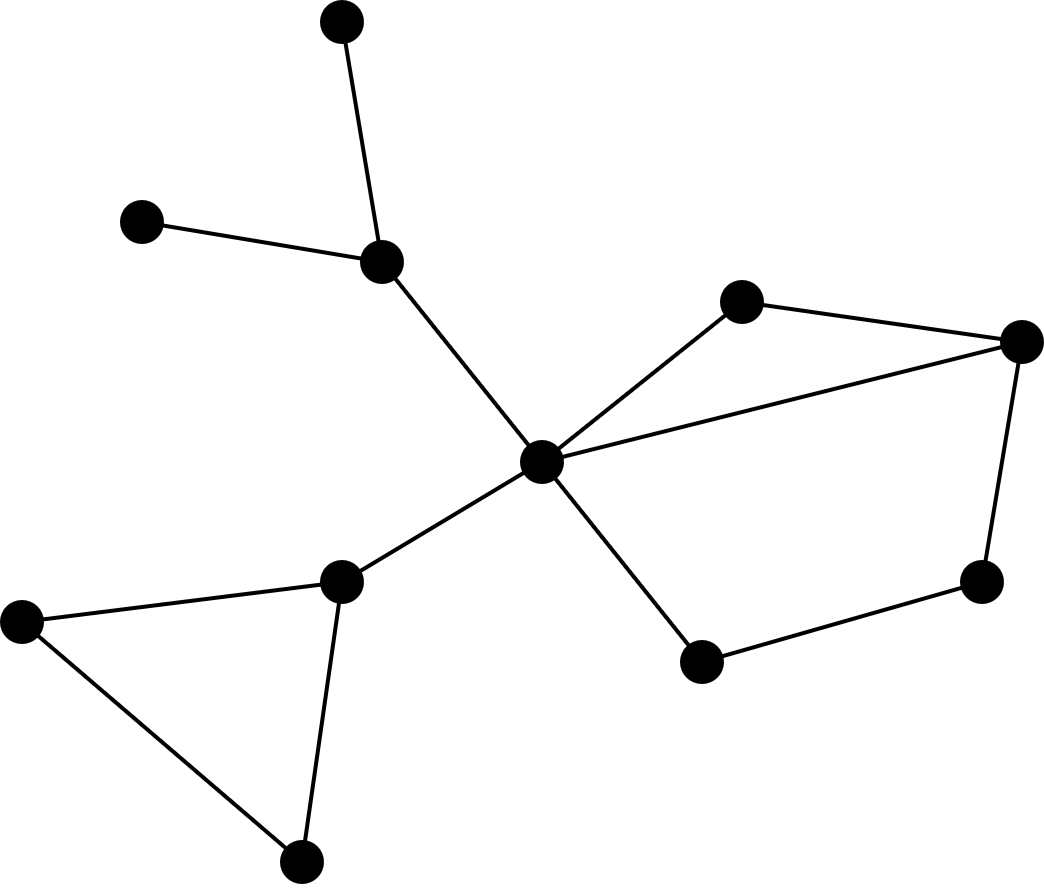}
            \captionof{figure}{An example graph \(H\)}
            \label{fig:H1}
        \end{minipage}%
        \begin{minipage}{.5\textwidth}
            \centering
            \includegraphics[width=.8\linewidth]{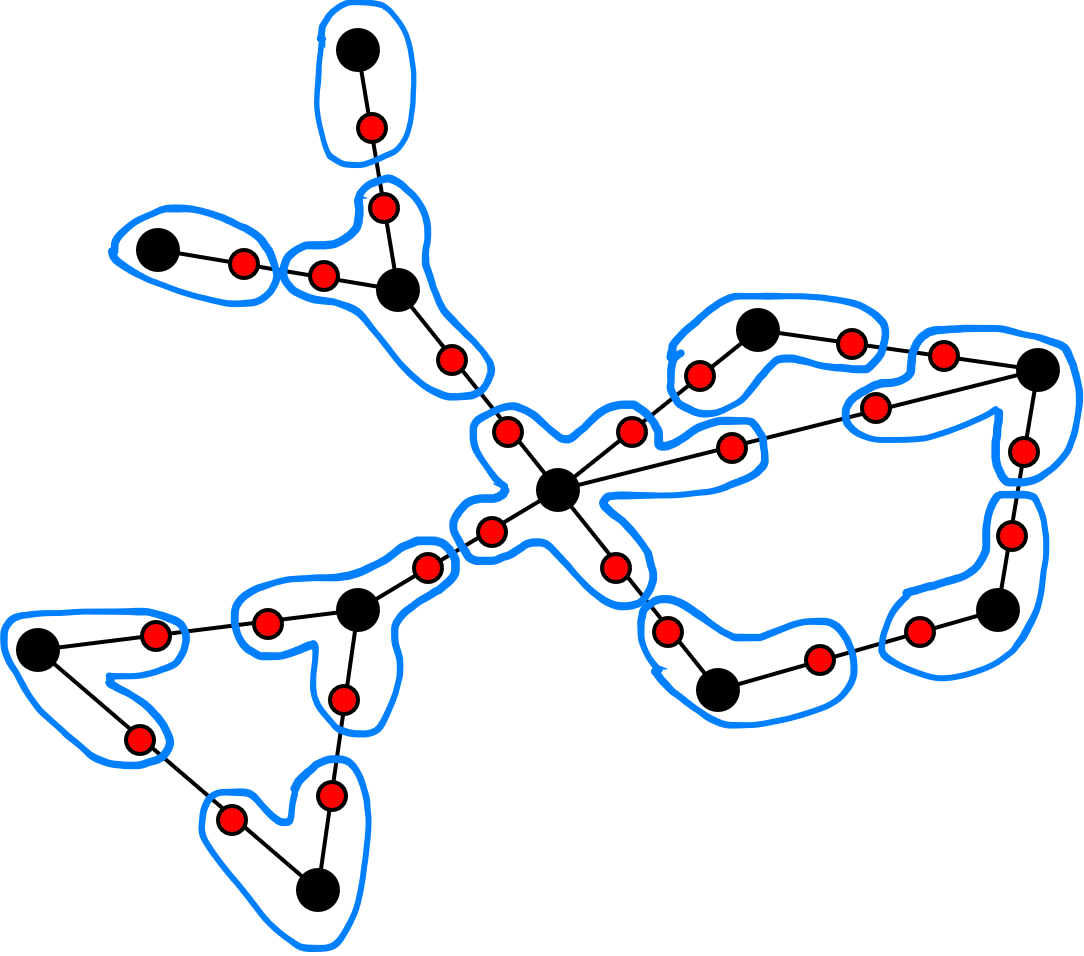}
            \captionof{figure}{The constructed graph \(G\)}
            \label{fig:G1}
        \end{minipage}
    \end{figure}

    Figures~\ref{fig:H1} and~\ref{fig:G1} show an example graph \(H\) containing a cycle and the corresponding graph \(G\). The red vertices represent the additional vertices inserted during subdivision, while the blue loops indicate the vertex subsets that induce \(H\) in \(G^\circ\).

    Notice that in the proof of Lemma~\ref{lem:podind1}, it would already be enough to subdivide every edge into a path of length \(2\). The construction used above is slightly easier to visualise.

    \begin{lemma}\label{lem:podind2}
        For every graph \(H\) that is a forest but not a subgraph of a path, there exists a graph \(G\) such that \(G\) is \(H\)-free, but \(G^\circ\) is not \(H\)-free.
    \end{lemma}

    \begin{proof}
        Assume that \(H\) is a forest but is not a subgraph of a path. Then \(H\)
        \begin{enumerate}
            \item contains no cycle;
            \item contains at least one vertex of degree at least \(3\).
        \end{enumerate}
        It follows that \(H\) contains \(S_3\), the star with three leaves, as an induced subgraph.

        Let \(G\) be obtained from \(H\) by replacing every vertex \(u\in V(H)\) with \(d(u)\geq 3\) by a clique \(K_{d(u)}\). We denote the clique replacing \(u\) by \(K_u\). The replacement is performed so that every edge originally incident with \(u\) leaves a different vertex of \(K_u\).

        We first show that \(G\) is \(H\)-free. Consider any vertex \(v\in V(G)\) with \(d(v)\geq 3\), and let \(v_1,v_2,v_3\in N_G(v)\) be pairwise distinct neighbours of \(v\). Every vertex in \(G\) of degree at least \(3\) belongs to some clique replacing a high-degree vertex of \(H\) and has exactly one neighbour outside that clique. Hence, without loss of generality, at least two of \(v_1,v_2,v_3\), say \(v_1\) and \(v_2\), lie in the same clique. Therefore \(v_1v_2\in E(G)\). Since this holds for every choice of \(v\) and its three neighbours, \(G\) is \(S_3\)-free, and consequently \(H\)-free.

        We now show that \(G^\circ\) is not \(H\)-free. Let \(X^\circ\) be the set of vertices of \(G^\circ\) representing the original vertices of \(H\): if a vertex \(w\in V(H)\) was not replaced by a clique, we take \(\{w\}\), while if \(w\) was replaced by a clique, we take \(V(K_w)\). It is immediate that two vertices in \(X^\circ\) are adjacent in \(G^\circ\) if and only if the corresponding original vertices of \(H\) were adjacent. Hence \(G^\circ[X^\circ]\simeq H\), so \(G^\circ\) contains \(H\) as an induced subgraph.

        Thus we have shown that no forest \(H\) that is not a subgraph of a path is IS-invariant.
    \end{proof}

    \begin{figure}[h]
        \centering
        \begin{minipage}{.5\textwidth}
            \centering
            \includegraphics[width=.8\linewidth]{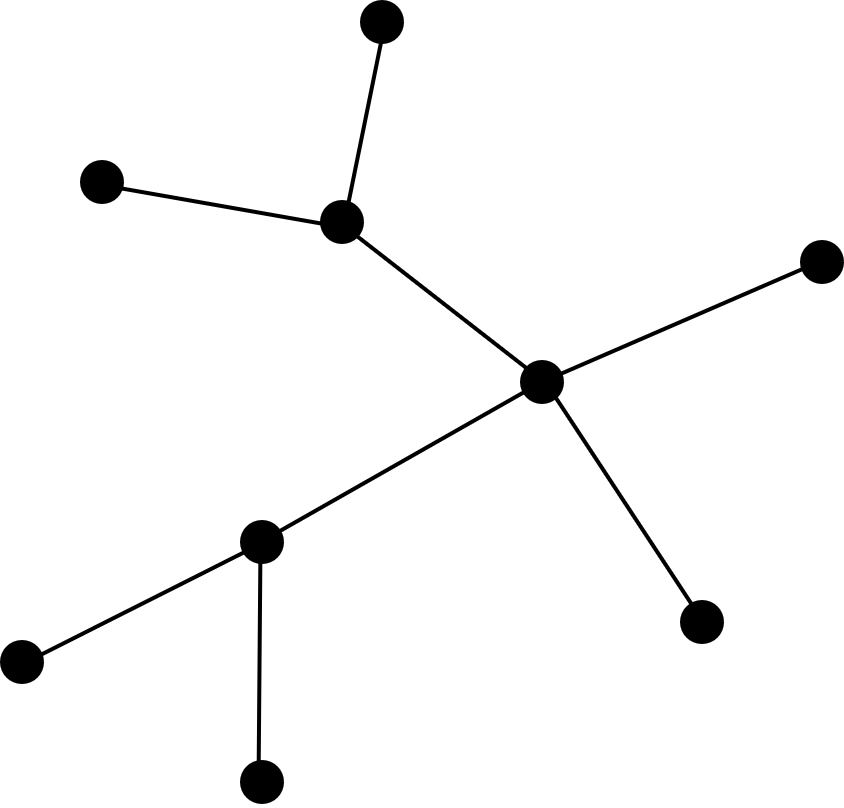}
            \captionof{figure}{An example graph \(H\)}
            \label{fig:H2}
        \end{minipage}%
        \begin{minipage}{.5\textwidth}
            \centering
            \includegraphics[width=.8\linewidth]{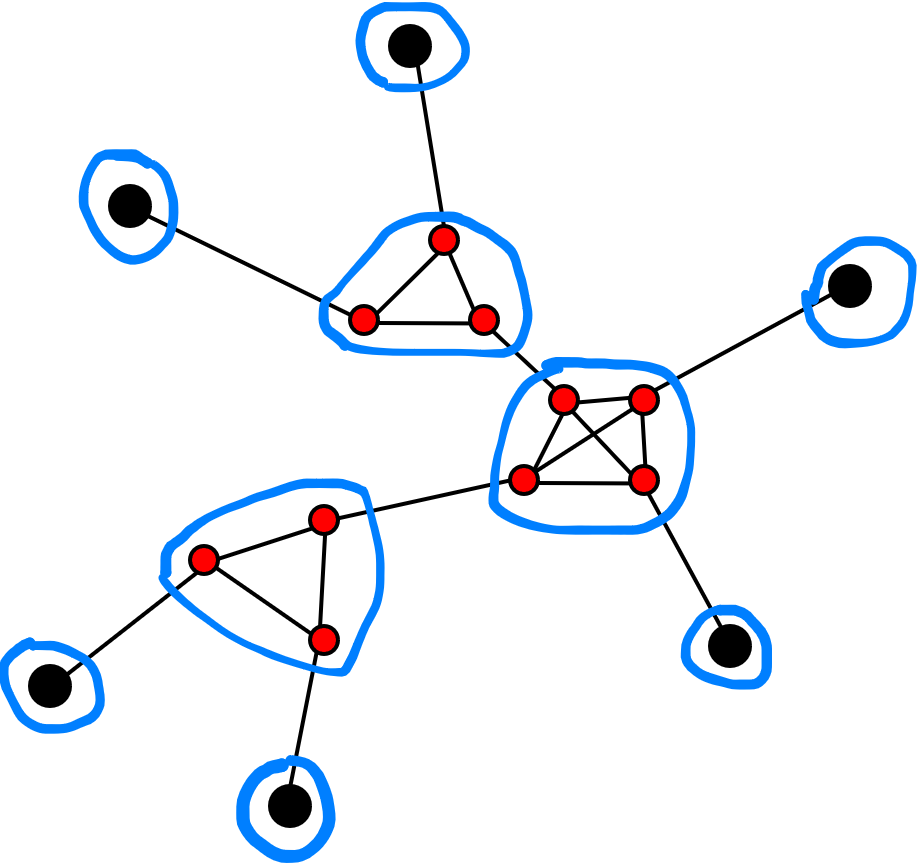}
            \captionof{figure}{The constructed graph \(G\)}
            \label{fig:G2}
        \end{minipage}
    \end{figure}

    Figures~\ref{fig:H2} and~\ref{fig:G2} show an example graph \(H\) that is a forest, in fact a tree, but is not a subgraph of a path, together with the corresponding graph \(G\). The red vertices represent vertices of the cliques into which the relevant vertices of \(H\) were split, while the blue loops indicate the vertex subsets that induce \(H\) in \(G^\circ\).

    It is worth noting that the construction from Lemma~\ref{lem:podind2} also gives counterexamples for a slightly larger class of graphs, namely for all graphs \(H\) containing \(S_3\) as an induced subgraph. This fact will be used later.

    By Lemmas~\ref{lem:podind1} and~\ref{lem:podind2}, for every graph \(H\) that is not a subgraph of a path, there exists a graph \(G\) such that \(G\) is \(H\)-free, but \(G^\circ\) is not \(H\)-free. Therefore, if a graph \(H\) is IS-invariant, then \(H\) must be a subgraph of a path. This proves the forward implication.
\end{proof}

\begin{proof}[Proof of the reverse implication]
    Let \(H\) be a subgraph of a path, and suppose that for some graph \(G\), the graph \(G^\circ\) contains \(H\) as an induced subgraph. We prove that \(G\) also contains \(H\) as an induced subgraph. In the rest of the proof, \(H'\) denotes an induced subgraph of \(G^\circ\) isomorphic to the given graph \(H\).

    Since \(H'\) is a subgraph of a path, we may write it as a disjoint union
    \[
        P_1,P_2,\dots,P_s,
    \]
    where every \(P_i\) is a path and is a maximal connected component of \(H'\).

    \begin{lemma}\label{lem:podind3}
        Let \(S,K\subseteq V(H')\), \(S\neq K\), induce distinct maximal connected components of \(H'\). Then, for every
        \[
            u\in V\!\left(G\left[\bigcup V(S)\right]\right)
            \quad\text{and}\quad
            v\in V\!\left(G\left[\bigcup V(K)\right]\right),
        \]
        there is no \(u\)-\(v\) edge in \(G\).
    \end{lemma}

    \begin{proof}
        Suppose otherwise. Let \(u\in V(G[\bigcup V(S)])\) and \(v\in V(G[\bigcup V(K)])\), and suppose that \(G\) contains an \(u\)-\(v\) edge. Then either the components \(S\) and \(K\) have non-empty intersection, or there is an edge of \(G\) between a vertex belonging to a blob from \(S\) and a vertex belonging to a blob from \(K\). By Definition~\ref{def:blob}, this gives an \(S\)-\(K\) edge in \(G^\circ\), contradicting the fact that \(S\) and \(K\) induce distinct maximal connected components of \(H'\).
    \end{proof}

    Thus, in \(G\) there is no edge between the parts corresponding to different paths \(P_i\). Therefore, it remains only to prove that inside each path \(P_1,\dots,P_s\), the graph \(G\) contains an induced path of at least the same length.

    \begin{lemma}\label{lem:podind4}
        Let \(P\coloneq \{V_1,V_2,\dots,V_n\}\) be an induced path in \(G^\circ\) on \(n\) vertices. Then \(G[\bigcup_{k=1}^n V_k]\) contains an induced path on at least \(n\) vertices.
    \end{lemma}

    \begin{proof}
        The graph \(G[\bigcup_{k=1}^n V_k]\) is connected, because for every \(i\in[n-1]\), the blobs \(V_i\) and \(V_{i+1}\) either intersect or are adjacent. Moreover, for every \(j\in\{2,3,\dots,n-1\}\), the graph \(G[\bigcup_{k=1,\ k\neq j}^n V_k]\) is not connected. Indeed, if it were connected, then there would be an edge between \(\bigcup_{k=1}^{j-1}V_k\) and \(\bigcup_{k=j+1}^{n}V_k\), which would contradict the fact that \(P\) is an induced path in \(G^\circ\).

        We shall use the following elementary fact.

        \begin{lemma}\label{lem:podind5}
            Let \(u\) and \(v\) be vertices of a connected graph \(F\), and let \(P_F\) be a shortest \(u\)-\(v\) path in \(F\). Then \(P_F\) is an induced path.
        \end{lemma}

        \begin{proof}
            Suppose not. Write \(P_F=\{u=v_1,v_2,\dots,v_m=v\}\). If there exist \(i,j\in[m]\) with \(i\leq j-2\) such that \(v_iv_j\in E(F)\), then \(\{v_1,\dots,v_i,v_j,\dots,v_m\}\) is a shorter \(u\)-\(v\) path in \(F\), a contradiction.
        \end{proof}

        Choose \(u\in V_1\) and \(v\in V_n\), and let
        \[
            P'=\{u=v_1,v_2,\dots,v_m=v\}
        \]
        be a shortest \(u\)-\(v\) path in \(G[\bigcup_{k=1}^n V_k]\). By Lemma~\ref{lem:podind5}, \(P'\) is induced. It remains to show that \(|P'|\geq n\).

        For \(j\in\{2,3,\dots,n-1\}\), let
        \[
            X_j \coloneq V_j-\bigcup_{k=1,\ k\neq j}^{n}V_k.
        \]
        Thus \(X_j\) consists of the vertices of \(G\) that, among the blobs on the path, belong only to \(V_j\).

        Suppose that \(X_j=\emptyset\) for some \(j\in\{2,3,\dots,n-1\}\). Then \(V_j\subseteq\bigcup_{k=1,\ k\neq j}^{n}V_k\), so
        \[
            G\!\left[\bigcup_{k=1,\ k\neq j}^{n}V_k\right]
            =
            G\!\left[\bigcup_{k=1}^{n}V_k\right].
        \]
        But the graph on the left is not connected, whereas the graph on the right is connected, a contradiction.

        Therefore, for every \(j\in\{2,3,\dots,n-1\}\), the path \(P'\) contains a vertex belonging to \(X_j\). The sets \(X_i\) and \(X_j\) are pairwise disjoint for \(i\neq j\), so \(P'\) contains at least \(n-2\) such vertices. Since the endpoints \(u\) and \(v\) do not belong to \(\bigcup_{j=2}^{n-1}X_j\), the path \(P'\) has at least \(n\) vertices.
    \end{proof}

    By Lemmas~\ref{lem:podind3} and~\ref{lem:podind4}, in \(G\) there are pairwise disjoint induced paths \(P'_1,P'_2,\dots,P'_s\) corresponding to the paths \(P_1,P_2,\dots,P_s\) in \(H'\), such that \(|P'_i|\geq |P_i|\). In particular, each \(P'_i\) contains an induced subpath on exactly \(|P_i|\) vertices. Therefore \(G\) contains \(H\) as an induced subgraph. This proves the reverse implication and completes the proof of Theorem~\ref{theo:podind}.
\end{proof}

\section{Induced topological minors}\label{mintop}

\begin{definition}[\textbf{ITM-invariant}]\label{def:mintop}
    A graph \(H\) is called \textbf{ITM-invariant} if for every graph \(G\), the fact that \(G\) does not contain \(H\) as an induced topological minor implies that \(G^\circ\) does not contain \(H\) as an induced topological minor.
\end{definition}

\begin{theorem}\label{theo:mintop}
    A graph \(H\) is ITM-invariant if and only if \(\Delta(H)\leq 2\) and \(H\) does not contain \(K_3\) as a subgraph.
\end{theorem}

\begin{proof}[Proof of the forward implication]
    Suppose that a graph \(H\) is ITM-invariant. Assume that \(\Delta(H)\geq 3\), or that \(H\) contains \(K_3\) as a subgraph.

    We first consider the case in which \(\Delta(H)\geq 3\) and \(H\) contains \(S_3\), the star with three leaves, as an induced subgraph.

    \begin{lemma}\label{lem:mintop1}
        For every graph \(H\) containing \(S_3\) as an induced subgraph, there exists a graph \(G\) such that \(G\) does not contain \(H\) as an induced topological minor, but \(G^\circ\) contains \(H\) as an induced topological minor.
    \end{lemma}

    \begin{proof}
        Of course, \(\Delta(H)\geq 3\), because \(S_3\) contains a vertex of degree \(3\), and the operations used to obtain an induced topological minor cannot create a vertex of degree larger than one already present in the graph.

        Construct \(G\) from \(H\) as in the proof of Lemma~\ref{lem:podind2}: replace every vertex \(u\in V(H)\) with \(d(u)\geq 3\) by a clique \(K_{d(u)}\), denoted \(K_u\), in such a way that every edge originally incident with \(u\) leaves a different vertex of \(K_u\).

        In the proof of Lemma~\ref{lem:podind2}, we showed that \(G^\circ\) contains \(H\) as an induced subgraph. In particular, \(G^\circ\) contains \(H\) as an induced topological minor. It remains to prove that \(G\) does not contain \(H\) as an induced topological minor.

        Let \(v\in V(G)\) be any vertex with \(d(v)\geq 3\), and let \(v_1,v_2,v_3\in N_G(v)\) be pairwise distinct neighbours of \(v\). Every vertex of degree at least \(3\) in \(G\) belongs to one of the cliques introduced above and has exactly one neighbour outside that clique. Hence, without loss of generality, \(v_1,v_2\in K_{v'}\) for the relevant clique \(K_{v'}\). Therefore \(v_1v_2\in E(G)\). Moreover, the involved vertices have degree at least \(3\), so the admissible operations for induced topological minors cannot remove this obstruction by simply suppressing one of these edges. It follows that \(G\) does not contain \(S_3\) as an induced topological minor, and hence does not contain \(H\) as an induced topological minor.

        Thus \(G\) demonstrates that no graph \(H\) containing \(S_3\) as an induced subgraph is ITM-invariant.
    \end{proof}

    Recall that an example graph \(H\) satisfying the assumptions of the lemma and the corresponding graph \(G\) are shown in Figures~\ref{fig:H2} and~\ref{fig:G2}.

    To finish the proof of the forward implication, it remains to rule out graphs \(H\) satisfying one of the following conditions:
    \begin{itemize}
        \item \(\Delta(H)\geq 3\) and \(H\) is \(S_3\)-free;
        \item \(H\) contains \(K_3\) as a subgraph.
    \end{itemize}

    Consider a graph \(H\) satisfying the first condition, that is, \(\Delta(H)\geq 3\) and \(H\) is \(S_3\)-free. Let \(v\in V(H)\) be a vertex with \(d(v)\geq 3\), and let \(v_1,v_2,v_3\in N_H(v)\) be pairwise distinct neighbours of \(v\). If none of the edges \(v_1v_2\), \(v_2v_3\), \(v_1v_3\) existed in \(H\), then \(H\) would contain an induced \(S_3\), a contradiction. Therefore, without loss of generality, \(v_1v_2\in E(H)\). Then the set \(\{v,v_1,v_2\}\) induces a triangle in \(H\). Hence every remaining case contains a triangle.

    \begin{lemma}\label{lem:mintop2}
        For every graph \(H\) containing \(K_3\) as a subgraph, there exists a graph \(G\) such that \(G\) does not contain \(H\) as an induced topological minor, but \(G^\circ\) contains \(H\) as an induced topological minor.
    \end{lemma}

    \begin{proof}
        Let \(T\) be a triangle in \(H\), with \(V(T)=\{v_1,v_2,v_3\}\).

        We construct \(G\) as follows. The vertex set \(V(G)\) contains all vertices in \(V(H)-V(T)\), together with four new vertices \(v'_1,v'_2,v'_3\) and \(s\). The edge set \(E(G)\) contains all edges between pairs of vertices of \(H\) that remain in \(G\), and additionally:
        \begin{itemize}
            \item for every \(i\in[3]\), all edges between \(v'_i\) and all vertices of \(N_H(v_i)\) that belong to \(G\);
            \item the three edges \(v'_1s\), \(v'_2s\), and \(v'_3s\).
        \end{itemize}

        Equivalently, \(G\) is obtained from \(H\) by replacing the triangle \(T\) with an induced star
        \[
            S,\qquad V(S)=\{v'_1,v'_2,v'_3,s\},
        \]
        with three leaves. Each leaf corresponds to one vertex of \(T\), and the centre \(s\) is adjacent only to the leaves.

        We first show that \(G\) does not contain \(H\) as an induced topological minor. We compare the numbers of vertices and edges in \(G\) and \(H\). We have \(|V(G)|=|V(H)|+1\), because the triangle on three vertices was replaced by a star on four vertices. Moreover, \(|E(G)|=|E(H)|\), because both \(T\) and \(S\) contain three edges, and all other edges are preserved.

        Since \(|V(G)|=|V(H)|+1\), obtaining \(H\) from \(G\), if possible, would require exactly one operation of the induced topological minor relation. Both vertex deletion and an admissible suppression or contraction decrease the number of edges by at least one. The exceptional case of deleting an isolated vertex is irrelevant here, because it cannot create \(H\) from \(G\). Therefore, for every graph \(G'\) obtained from \(G\) by one such operation, we have \(|E(G')|<|E(H)|\). Hence \(G'\neq H\), so \(G\) does not contain \(H\) as an induced topological minor.

        We now show that \(G^\circ\) contains \(H\) as an induced topological minor. Consider
        \[
            V^\circ
            =
            \left(\bigcup_{v\in V(H)-V(T)}\{\{v\}\}\right)
            \cup
            \bigl\{\{v'_1,s\},\{v'_2,s\},\{v'_3,s\}\bigr\}.
        \]
        The vertices \(\{v'_1,s\}\), \(\{v'_2,s\}\), and \(\{v'_3,s\}\) induce a triangle in \(G^\circ\). Moreover, for every \(i\in[3]\), the set \(\{v'_i,s\}\) is adjacent exactly to those singleton vertices in \(V^\circ\) corresponding to vertices adjacent to \(v_i\) in \(H\). Hence \(G^\circ[V^\circ]\simeq H\). Thus \(G^\circ\) contains \(H\) as an induced subgraph, and in particular as an induced topological minor.

        This demonstrates that no graph \(H\) containing \(K_3\) as a subgraph is ITM-invariant.
    \end{proof}

    \begin{figure}[h]
        \centering
        \begin{minipage}{.5\textwidth}
            \centering
            \includegraphics[width=.8\linewidth]{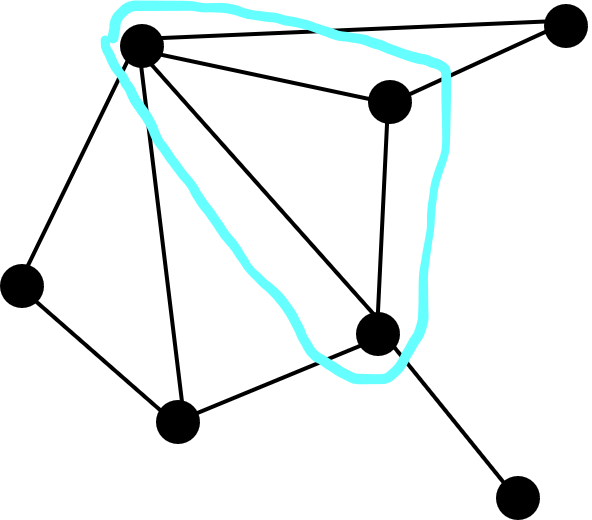}
            \captionof{figure}{An example graph \(H\)}
            \label{fig:H3}
        \end{minipage}%
        \begin{minipage}{.5\textwidth}
            \centering
            \includegraphics[width=.8\linewidth]{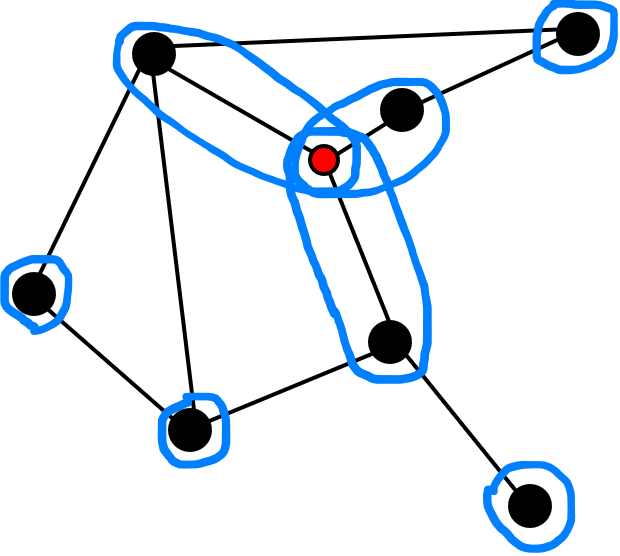}
            \captionof{figure}{The constructed graph \(G\)}
            \label{fig:G3}
        \end{minipage}
    \end{figure}

    Figures~\ref{fig:H3} and~\ref{fig:G3} show an example graph \(H\) containing a triangle, in particular the triangle \(T\) chosen in the construction, and the corresponding graph \(G\). The red vertex is the centre of the newly created star, while the blue loops indicate the vertex subsets that induce \(H\) in \(G^\circ\).

    Notice that the construction from Lemma~\ref{lem:mintop2} also gives a counterexample for the analogous statement concerning induced minors, because every induced-minor operation, except for special irrelevant cases, also decreases the number of edges.

    By Lemmas~\ref{lem:mintop1} and~\ref{lem:mintop2}, for every graph \(H\) such that \(\Delta(H)\geq 3\) or \(H\) contains \(K_3\) as a subgraph, there exists a graph \(G\) such that \(G\) does not contain \(H\) as an induced topological minor, while \(G^\circ\) contains \(H\) as an induced topological minor. Hence, if a graph \(H\) is ITM-invariant, then \(H\) has maximum degree at most \(2\) and contains no triangle. This proves the forward implication.
\end{proof}

\begin{proof}[Proof of the reverse implication]
    Let \(H\) be a graph with \(\Delta(H)\leq 2\), and assume that \(H\) contains no \(K_3\) as a subgraph. Suppose that, for some graph \(G\), the graph \(G^\circ\) contains \(H\) as an induced topological minor. We prove that \(G\) also contains \(H\) as an induced topological minor.

    Every maximal connected component of \(H\) is a path or a cycle, because \(\Delta(H)\leq 2\). Since \(H\) contains no triangle, every cyclic component has length at least \(4\). Thus we can write \(H\) as a disjoint union of paths and cycles
    \[
        P_1,\dots,P_s,C_1,\dots,C_k,
    \]
    where each \(P_i\) is a path component and each \(C_i\) is a cycle component of length at least \(4\).

    Since \(G^\circ\) contains \(H\) as an induced topological minor, it contains some subdivision \(H'\) of \(H\) as an induced subgraph. We can write \(H'\) as a disjoint union
    \[
        P'_1,\dots,P'_s,C'_1,\dots,C'_k,
    \]
    where \(P'_i\) is a subdivision of \(P_i\), and \(C'_j\) is a subdivision of \(C_j\). Each of these subgraphs induces a maximal connected component of \(H'\), and naturally \(|P'_i|\geq |P_i|\) and \(|C'_j|\geq |C_j|\).

    By Lemma~\ref{lem:podind3}, it is enough to show that each path component \(P'_i\) lifts to an induced path in \(G\) of at least the same length, and each cycle component \(C'_j\) lifts to an induced cycle in \(G\) of at least the same length. The path case follows from Lemma~\ref{lem:podind4}. It remains to handle cycles.

    \begin{lemma}\label{lem:mintop3}
        Let \(C\coloneq\{V_1,V_2,\dots,V_n\}\) be an induced cycle in \(G^\circ\) of length \(n\geq 4\). Then \(G[\bigcup_{k=1}^{n}V_k]\) contains an induced cycle of length at least \(n\).
    \end{lemma}

    \begin{proof}
        For every \(j\in[n]\), define
        \[
            X_j\coloneq V_j-\bigcup_{k=1,\ k\neq j}^{n}V_k.
        \]
        As in the proof of Lemma~\ref{lem:podind4}, every \(X_j\) is non-empty. Indeed, in an induced cycle of length at least \(4\), every vertex has two mutually non-adjacent neighbours; hence each blob \(V_j\) must contain at least one vertex separating the vertices belonging to the two neighbouring blobs.

        Moreover, since \(G[\bigcup V(C)]\) is connected, each \(X_j\) contains two non-empty subsets \(X_j^L\) and \(X_j^R\) such that there is an edge in \(G\) from \(X_j^L\) to \(V_{j-1}\) and an edge from \(X_j^R\) to \(V_{j+1}\), with indices taken cyclically. In particular, for \(X_n\) we choose \(X_n^L\) and \(X_n^R\) so that there are edges from \(X_n^L\) to \(V_{n-1}\) and from \(X_n^R\) to \(V_1\).

        Let \(P\) be a shortest path in \(G[X_n]\) that starts in \(X_n^L\) and ends in \(X_n^R\). By Lemma~\ref{lem:podind5}, \(P\) is induced. By minimality, \(P\) contains no vertices of \(X_n^L\) or \(X_n^R\) except for its endpoints, which we denote by \(u_p\) and \(v_p\).

        Now consider the path \(Q'=\{V_1,V_2,\dots,V_{n-1}\}\) in \(G^\circ\). Let
        \[
            V'_1\coloneq V_1\cap N(v_p),
            \qquad
            V'_{n-1}\coloneq V_{n-1}\cap N(u_p).
        \]
        Let \(Q\) be a shortest path in \(G[\bigcup_{k=1}^{n-1}V_k]\) with one endpoint in \(V'_1\) and the other in \(V'_{n-1}\). By Lemma~\ref{lem:podind4}, \(Q\) is an induced path of length at least \(n-1\). Moreover, by minimality of \(Q\), it contains no vertices of \(V'_1\) or \(V'_{n-1}\) except for its endpoints, which we denote by \(u_q\) and \(v_q\).

        The paths \(P\) and \(Q\) are induced paths in the corresponding induced subgraphs of \(G\). In \(G[\bigcup_{k=1}^{n}V_k]\), the only edges between \(P\) and \(Q\) are \(v_pu_q\) and \(v_qu_p\). Since \(|Q|\geq n-1\) and \(|P|\geq 1\), the graph \(G[P\cup Q]\) is an induced cycle of length at least \(n\).
    \end{proof}

    We have now shown that for every maximal connected component of \(H'\), if it is a path \(P'_i\), then \(G\) contains an induced path \(P''_i\) with \(|P''_i|\geq |P'_i|\geq |P_i|\), and if it is a cycle \(C'_i\), then \(G\) contains an induced cycle \(C''_i\) with \(|C''_i|\geq |C'_i|\geq |C_i|\). Therefore \(G\) contains \(H\) as an induced topological minor. This proves the reverse implication and completes the proof of Theorem~\ref{theo:mintop}.
\end{proof}

\section{Induced minors}\label{minind}

\begin{definition}[\textbf{IM-invariant}]\label{def:minind}
    A graph \(H\) is called \textbf{IM-invariant} if for every graph \(G\), the fact that \(G\) does not contain \(H\) as an induced minor implies that \(G^\circ\) does not contain \(H\) as an induced minor.
\end{definition}

\begin{theorem}\label{theo:minind}
    A graph \(H\) is IM-invariant if and only if \(H\) satisfies both of the following conditions:
    \begin{itemize}
        \item \(H\) contains no \(K_3\), that is, no triangle, as a subgraph;
        \item for every \(v\in V(H)\) with \(d(v)\geq 3\), all neighbours of \(v\) in \(H\) have degree exactly \(2\).
    \end{itemize}
\end{theorem}

\begin{proof}[Proof of the forward implication]
    Suppose that a graph \(H\) is IM-invariant. Assume that \(H\) contains \(K_3\) as a subgraph, or that \(H\) contains a vertex \(v\) with \(d(v)\geq 3\) that has a neighbour of degree different from \(2\).

    We first consider the case in which \(H\) contains \(K_3\) as a subgraph.

    \begin{lemma}\label{lem:minind1}
        For every graph \(H\) containing a triangle, there exists a graph \(G\) such that \(G\) does not contain \(H\) as an induced minor, but \(G^\circ\) contains \(H\) as an induced minor.
    \end{lemma}

    \begin{proof}
        Construct \(G\) from \(H\) as in the proof of Lemma~\ref{lem:mintop2}: replace a chosen triangle \(K_3\) in \(H\) by a star \(S_3\), in such a way that every edge originally leaving the clique \(K_3\) leaves a different leaf of the star.

        In Lemma~\ref{lem:mintop2}, we showed that \(G^\circ\) contains \(H\) as an induced topological minor, and hence in particular as an induced minor. It remains only to show that \(G\) does not contain \(H\) as an induced minor.

        This part is analogous to the proof of Lemma~\ref{lem:mintop2}. We have \(|V(G)|=|V(H)|+1\) and \(|E(G)|=|E(H)|\). Obtaining \(H\) from \(G\) would require exactly one induced-minor operation. However, every such operation, except for deleting an isolated vertex, decreases the number of edges in \(G\). Deleting an isolated vertex cannot create \(H\) from this construction. Hence \(G\) does not contain \(H\) as an induced minor.

        Therefore \(G\) demonstrates that no graph \(H\) containing \(K_3\) as a subgraph is IM-invariant.
    \end{proof}

    Recall that an example graph \(H\) satisfying the assumptions of the lemma and the corresponding graph \(G\) are shown in Figures~\ref{fig:H3} and~\ref{fig:G3}.

    Before giving counterexamples for the remaining graphs \(H\), we prove the following invariant.

    \begin{lemma}\label{lem:minind2}
        No induced-minor operation on a graph \(G\) can increase the sum of the degrees of the vertices \(v\in V(G)\) with \(d(v)\geq 3\).
    \end{lemma}

    \begin{proof}
        We consider the two allowed operations.
        \begin{itemize}
            \item \textbf{Vertex deletion.} Deleting a vertex removes all edges incident with it. It does not increase the degrees of any other vertices. Therefore the sum of the degrees of vertices satisfying \(d(v)\geq 3\) can only remain unchanged or decrease.

            \item \textbf{Edge contraction.} When two adjacent vertices \(v_1,v_2\in V(G)\) are contracted into a new vertex \(v_n\), we have
                  \[
                      d(v_n)\leq d(v_1)+d(v_2)-2.
                  \]
                  This is because the edge \(v_1v_2\) disappears, and common neighbours of \(v_1\) and \(v_2\) become single neighbours of \(v_n\). If both \(v_1\) and \(v_2\) have degree smaller than \(3\), then \(v_n\) also has degree smaller than \(3\). If at least one of them, say \(v_1\), has \(d(v_1)\geq 3\), then either \(d(v_2)\geq 3\), in which case \(d(v_n)<d(v_1)+d(v_2)\), or \(d(v_2)\leq 2\), in which case \(d(v_n)\leq d(v_1)\).
        \end{itemize}

        In both cases, the sum of the degrees of vertices with degree at least \(3\) does not increase.
    \end{proof}

    We are now ready to show that a graph \(H\) cannot be IM-invariant whenever \(H\) contains a vertex \(v\) with \(d(v)\geq 3\) that has a neighbour whose degree is not \(2\).

    \begin{lemma}\label{lem:minind3}
        Let \(H\) be a graph containing a vertex \(v\) with \(d(v)\geq 3\) and a neighbour \(u\in N_H(v)\) such that \(d(u)\neq 2\). Then there exists a graph \(G\) such that \(G\) does not contain \(H\) as an induced minor, but \(G^\circ\) contains \(H\) as an induced minor.
    \end{lemma}

    \begin{proof}
        We construct the graph \(G\). Let \(v\) be the chosen vertex with \(d_H(v)=n+1\geq 3\), and let \(u\) be a neighbour of \(v\) such that \(d_H(u)=m+1\neq 2\).

        The vertex set of \(G\) consists of all vertices \(V(H)-\{u,v\}\), together with new vertices
        \[
            \{s,v_1,v_2,\dots,v_{n+m}\}.
        \]
        The edge set of \(G\) contains all edges between pairs of vertices of \(H\) that remain in \(G\), and additionally:
        \begin{itemize}
            \item for every \(i\in[n+m]\), the edge \(v_is\);
            \item for \(i\in[n]\), edges joining consecutive vertices \(v_i\) to distinct neighbours of \(v\) in \(V(H)-\{u\}\);
            \item for \(i\in[n+m]\setminus[n]\), edges joining consecutive vertices \(v_i\) to distinct neighbours of \(u\) in \(V(H)-\{v\}\).
        \end{itemize}

        Informally, \(G\) is obtained from \(H\) by replacing the chosen pair of vertices \(u,v\) by a star whose number of leaves is equal to the total number of neighbours of \(u\) and \(v\), excluding the edge \(uv\). Each leaf of this star is adjacent to a different former neighbour of either \(u\) or \(v\).

        We first show that \(G\) does not contain \(H\) as an induced minor. If \(d_H(u)\geq 3\), then in \(G\) the number of vertices of degree at least \(3\) decreases by \(1\): the two vertices \(u\) and \(v\) are replaced by one vertex \(s\). Moreover,
        \[
            d_G(s)=d_H(u)+d_H(v)-2.
        \]
        If \(d_H(u)=1\), then the number of vertices of degree at least \(3\) does not change, but \(d_G(s)=d_H(v)-1\). Thus, in both cases, the sum of the degrees of vertices of degree at least \(3\) is smaller in \(G\) than it was in \(H\). By Lemma~\ref{lem:minind2}, \(G\) cannot contain \(H\) as an induced minor.

        We now show that \(G^\circ\) contains \(H\) as an induced minor. In fact, it contains \(H\) as an induced subgraph. Construct a set \(V^\circ\subseteq V(G^\circ)\) consisting of:
        \begin{itemize}
            \item singleton sets containing one vertex of \(H\) that also belongs to \(G\);
            \item the two sets
                  \[
                      \{s\}\cup\{v_1,\dots,v_n\}
                      \quad\text{and}\quad
                      \{s\}\cup\{v_{n+1},\dots,v_{n+m}\}.
                  \]
        \end{itemize}
        The singleton vertices represent the unchanged vertices of \(H\). The two larger blobs represent \(v\) and \(u\), respectively. They are adjacent to each other, and they are adjacent exactly to the singleton vertices representing the corresponding neighbours of \(v\) and \(u\) in \(H\). Therefore \(G^\circ\) contains \(H\) as an induced subgraph, and hence as an induced minor. If \(d_H(u)=1\), then the second larger set is simply \(\{s\}\).

        This demonstrates that no graph \(H\) containing a vertex \(v\) of degree at least \(3\) with a neighbour of degree different from \(2\) is IM-invariant.
    \end{proof}

    \begin{figure}[h]
        \centering
        \begin{minipage}{.5\textwidth}
            \centering
            \includegraphics[width=.8\linewidth]{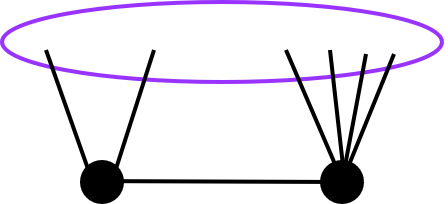}
            \captionof{figure}{An example graph \(H\)}
            \label{fig:H4}
        \end{minipage}%
        \begin{minipage}{.5\textwidth}
            \centering
            \includegraphics[width=.8\linewidth]{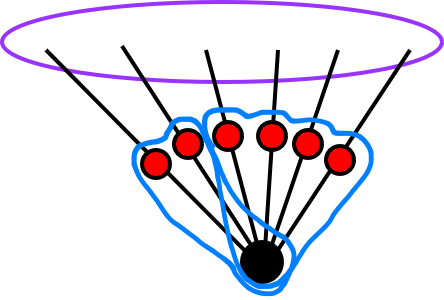}
            \captionof{figure}{The constructed graph \(G\)}
            \label{fig:G4}
        \end{minipage}
    \end{figure}

    \begin{figure}[h]
        \centering
        \begin{minipage}{.5\textwidth}
            \centering
            \includegraphics[width=.8\linewidth]{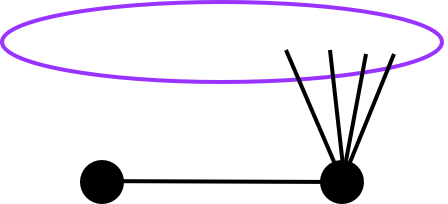}
            \captionof{figure}{An example graph \(H\)}
            \label{fig:H5}
        \end{minipage}%
        \begin{minipage}{.5\textwidth}
            \centering
            \includegraphics[width=.8\linewidth]{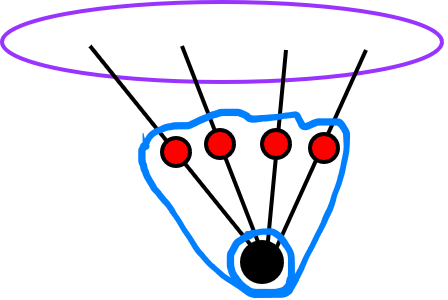}
            \captionof{figure}{The constructed graph \(G\)}
            \label{fig:G5}
        \end{minipage}
    \end{figure}

    Figures~\ref{fig:H4} and~\ref{fig:G4} show an example graph \(H\) containing a vertex \(v\) of degree at least \(3\) with a neighbour of degree different from \(2\), together with the constructed graph \(G\). The red vertices denote the newly created vertices \(\{v_1,v_2,\dots,v_{n+m}\}\), the purple ellipses denote the unchanged rest of the graph, and the blue loops denote the vertex subsets that induce \(H\) in \(G^\circ\), except for the unchanged rest of the graph where the natural mapping is used. Figures~\ref{fig:H5} and~\ref{fig:G5} show the special case in which the chosen vertex \(v\) of degree at least \(3\) has a neighbour of degree \(1\).

    By Lemmas~\ref{lem:minind1} and~\ref{lem:minind3}, for every graph \(H\) such that \(H\) contains \(K_3\) as a subgraph, or \(H\) has a vertex of degree at least \(3\) with a neighbour of degree different from \(2\), there exists a graph \(G\) such that \(G\) does not contain \(H\) as an induced minor, but \(G^\circ\) contains \(H\) as an induced minor. Therefore, if a graph \(H\) is IM-invariant, then \(H\) contains no triangle and every neighbour of every vertex of degree at least \(3\) has degree exactly \(2\). This proves the forward implication.
\end{proof}

Before proving the reverse implication of Theorem~\ref{theo:minind}, we establish the following useful property.

\begin{theorem}\label{theo:min-pod}
    Let \(G^\circ\) be the blob graph of a graph \(G\), and let \(H\) be a graph containing no \(K_3\) as a subgraph. If \(G^\circ\) contains \(H\) as an induced minor, then \(G^\circ\) contains \(H\) as an induced subgraph.
\end{theorem}

\begin{proof}
    Let \(H\) be an induced minor of \(G^\circ\). Then every vertex of \(H\) is obtained by contracting some connected set of vertices of \(G^\circ\). Write the vertices of \(H\) as
    \[
        V(H)=
        \bigl\{
        V_1\coloneq (V_{1,1},V_{1,2},\dots,V_{1,k_1}),
        \ldots,
        V_m\coloneq (V_{m,1},V_{m,2},\dots,V_{m,k_m})
        \bigr\},
    \]
    where \(V_i\) is obtained by contracting the vertices \(V_{i,1},V_{i,2},\dots,V_{i,k_i}\) of \(G^\circ\). Of course, for \((i_1,j_1)\neq(i_2,j_2)\), we have \(V_{i_1,j_1}\neq V_{i_2,j_2}\), because every vertex of \(G^\circ\) may be contracted into at most one vertex of \(H\).

    Observe that two vertices \(V_i,V_j\in V(H)\) are adjacent if and only if the sets
    \[
        \bigcup_{s=1}^{k_i}V_{i,s}
        \quad\text{and}\quad
        \bigcup_{t=1}^{k_j}V_{j,t}
    \]
    are adjacent as vertices of \(G^\circ\), or are equal. Indeed, \(V_i\) and \(V_j\) are adjacent in \(H\) precisely when there is an edge in \(G^\circ\) between some \(V_{i,s}\) and \(V_{j,t}\). This means exactly that \(V_{i,s}\cap V_{j,t}\neq\emptyset\), or that there is an edge in \(G\) between the corresponding subsets. These are precisely the adjacency conditions for the two unions in \(G^\circ\), with equality treated as the exceptional case.

    For each \(i\in[m]\), put
    \[
        X_i\coloneq \bigcup_{s=1}^{k_i}V_{i,s}.
    \]
    Then, in \(G^\circ\), the vertices \(X_i\) and \(X_j\) are adjacent if and only if the corresponding vertices \(V_i\) and \(V_j\) are adjacent in \(H\), except in the special case \(X_i=X_j\).

    Suppose that \(X_i=X_j\) for some \(i\neq j\). Then \(V_i\) and \(V_j\) are adjacent in \(H\). If one of them, say \(V_i\), had a neighbour \(V_k\) different from \(V_j\), then \(V_j\) would also be adjacent to \(V_k\), which would create a triangle in \(H\), a contradiction. Similarly, no third vertex \(V_k\) can satisfy \(X_k=X_i\). Hence \(V_i\) and \(V_j\) form a connected component \(K_2\) of \(H\), disjoint from the rest of the graph.

    Moreover, \(|X_i|=|X_j|\geq 2\). If this set had size \(1\), it could not arise in two different ways as the above union while preserving the contracted connected sets. Since \(G[X_i]\) is connected and has at least two vertices, we may choose a vertex \(x\in X_i\) such that \(X_i-\{x\}\) is still connected; for example, take a non-cut vertex of the connected graph \(G[X_i]\).

    Let \(X\coloneq\bigcup_{i=1}^{m}X_i\). For every pair \(i,j\) with \(i\neq j\) and \(X_i=X_j\), add to \(X\) the vertex of \(G^\circ\)
    \[
        X'_i\coloneq X_i-\{x\},
    \]
    where \(x\) is chosen as above. The only vertex among the original \(X_\ell\)'s adjacent to \(X'_i\) is \(X_i\). Furthermore, \(X_i\) has no neighbours among the remaining \(X_\ell\)'s, because otherwise \(H\) would contain a triangle, as shown above.

    Therefore the induced subgraph of \(G^\circ\) on the resulting set of vertices is isomorphic to \(H\). The ordinary edges of \(H\) are represented by the corresponding adjacencies between the \(X_i\)'s, while each isolated \(K_2\) component created by a pair with the same union is represented by \(X_i\) and \(X'_i\).
\end{proof}

\begin{proof}[Proof of the reverse implication of Theorem~\ref{theo:minind}]
    Let \(H\) be a graph such that \(H\) contains no \(K_3\) as a subgraph and, for every vertex \(v\in V(H)\) with \(d(v)\geq 3\), all neighbours of \(v\) have degree exactly \(2\). Suppose that, for some graph \(G\), the blob graph \(G^\circ\) contains \(H\) as an induced minor. We prove that \(G\) also contains \(H\) as an induced minor.

    Since \(H\) is triangle-free, Theorem~\ref{theo:min-pod} implies that \(G^\circ\) contains \(H\) as an induced subgraph. Let \(H'\) be an induced subgraph of \(G^\circ\) isomorphic to \(H\), and let \(S'_1,S'_2,\dots,S'_k\) be all maximal connected components of \(H'\).

    By Lemma~\ref{lem:podind3}, it is enough to prove, for every component \(S'_i\), that \(G\) contains a copy of the corresponding component as an induced minor. If \(\Delta(S'_i)\leq 2\), then \(S'_i\) is a path or a cycle. If it is a cycle, it has length at least \(4\), because \(H'\) is triangle-free. By Lemmas~\ref{lem:podind4} and~\ref{lem:mintop3}, respectively, \(G\) contains an induced path or cycle of length at least \(|S'_i|\) in the vertices belonging to \(S'_i\). By contracting along that induced path or cycle as needed, we obtain \(S'_i\) as an induced minor. Therefore it remains to handle components \(S'_i\) with \(\Delta(S'_i)\geq 3\).

    \begin{lemma}\label{lem:minind5}
        Let \(S\) be a graph satisfying:
        \begin{itemize}
            \item \(S\) is connected;
            \item \(\Delta(S)\geq 3\);
            \item \(S\) contains no \(K_3\) as a subgraph;
            \item for every \(v\in V(S)\) with \(d(v)\geq 3\), all neighbours of \(v\) in \(S\) have degree exactly \(2\).
        \end{itemize}
        If \(G^\circ\) contains \(S\) as an induced subgraph, then
        \[
            G\!\left[\bigcup_{X\in V(S)}X\right]
        \]
        contains \(S\) as an induced minor.
    \end{lemma}

    \begin{proof}
        Let \(T\subseteq V(S)\) be the set of vertices of \(S\) whose degree is different from \(2\). Let
        \[
            \mathcal{P}\coloneq\{P_1,P_2,\dots,P_l\}
        \]
        be the set of all maximal paths in \(S\) whose vertices do not belong to \(T\). Equivalently, \(\mathcal{P}\) contains all maximal paths in \(S\) whose every vertex \(V\) satisfies \(d_S(V)=2\).

        At this point, observe that every path in \(\mathcal{P}\) is an induced path in \(S\). This follows from the fact that every internal vertex of such a path has degree \(2\) both in the path and in \(S\), while the endpoints cannot be adjacent; otherwise \(S\) would be a cycle, contradicting \(\Delta(S)\geq 3\).

        Moreover, every vertex \(V\in V(S)\) with \(d_S(V)=2\) lies on some path in \(\mathcal{P}\). If not, then starting a depth-first search from \(V\), we would never encounter a vertex of degree larger than \(2\). Then \(S\) would be a cycle, again contradicting \(\Delta(S)\geq 3\).

        Take any path \(P_i\in\mathcal{P}\), say
        \[
            P_i=\{V_1,V_2,\dots,V_n\}.
        \]
        In \(S\), there exist vertices \(U,V\in T\) such that \(V_1\in N(U)\) and \(V_n\in N(V)\). Otherwise \(P_i\) would not be maximal, or \(S\) would be a cycle. It is possible that \(U=V\).

        We have \(V_1-U\neq\emptyset\) and \(V_n-V\neq\emptyset\). If this failed, say for \(U\), then \(U\) would be adjacent to \(V_2\), or to \(V\) when \(|P_i|=1\), and \(S\) would contain a triangle.

        Define a new sequence
        \[
            P'_i=\{V'_1,V'_2,\dots,V'_n\}
        \]
        by setting \(V'_j\coloneq V_j\) for \(2\leq j\leq n-1\). We select \(V'_1\) to be a maximal subset of \(V_1-U\) inducing a complete subgraph such that \(V'_1\in N_{G^\circ}(V_2)\). Symmetrically, let \(V'_n\) be a maximal subset of \(V_n-V\) inducing a complete subgraph such that \(V'_n\in N_{G^\circ}(V_{n-1})\). The sequence \(P'_i\) remains an induced path in \(G^\circ\).

        For every \(j\in[n]\), let
        \[
            X_j\coloneq V'_j-\bigcup_{k=1,\ k\neq j}^{n}V'_k.
        \]
        As in the proof of Lemma~\ref{lem:podind4}, for \(j\in\{2,\dots,n-1\}\) we have \(X_j\neq\emptyset\), because every vertex of \(P'_i\) has two mutually non-adjacent neighbours. Furthermore, \(X_1\neq\emptyset\) and \(X_n\neq\emptyset\), because in \(S\) there must be vertices separating \(\{V'_2,V'_3,\dots,V'_{n-1}\}\) from \(U\) and \(V\), respectively.

        Consider a shortest path in \(G[\bigcup_{j=1}^{n}V'_j]\) that starts in \(X_1\cap N(U)\) and ends in \(X_n\cap N(V)\). These two sets are non-empty; otherwise \(S\) would not contain the edges \(UV_1\) and \(VV_n\). This path is a shortest path between some \(u\in V'_1\) and \(v\in V'_n\), so by Lemma~\ref{lem:podind5} it is induced. It must contain at least one vertex from each \(X_j\), \(j\in[n]\). Thus it is an induced path on at least \(n\) vertices. For the given \(P_i\), denote this constructed path in \(G[\bigcup_{X\in V(S)}X]\) by \(p_i\).

        We are ready to show that \(G[\bigcup_{X\in V(S)}X]\) contains \(S\) as an induced minor. First, using vertex deletions, keep only vertices that either belong to \(\bigcup_{Y\in T}Y\) or to one of the paths \(p_i\). Each \(p_i\) is already an induced path of length at least \(|P_i|\). Hence, using contractions, we may shorten every \(p_i\) to a path \(p'_i\) of exactly the same length as \(P_i\). Next, for each \(Y\in T\), contract all vertices of \(Y\) into one vertex.

        The resulting graph has one vertex representing every \(Y\in T\), adjacent to the appropriate number of induced paths of the same lengths as the paths adjacent to \(Y\) in \(S\). Every degree-\(2\) vertex of \(S\) lies on one of the paths \(P_i\), and hence is represented on the corresponding path \(p'_i\). Therefore \(G[\bigcup_{X\in V(S)}X]\) contains \(S\) as an induced minor.
    \end{proof}

    We have shown that for every maximal connected component of \(H'\), the graph \(G\) contains a copy of that component as an induced minor. By Lemma~\ref{lem:podind3}, the components do not interact, so \(G\) contains \(H\) as an induced minor. This proves the reverse implication and completes the proof of Theorem~\ref{theo:minind}.
\end{proof}

\section{Conclusion}

We have determined for which graphs \(H\), when \(H\) is forbidden as one of the considered substructures, the corresponding blob graph also avoids \(H\) as the same substructure. We obtained characterisations for every substructure considered in the project, thereby achieving the goals of the project.

\bibliographystyle{plain}

\bibliography{bibliography.bib}

@misc{gartland2021findinglargeinducedsparse,
      title={Finding large induced sparse subgraphs in $C_{>t}$-free graphs in quasipolynomial time}, 
      author={Peter Gartland and Daniel Lokshtanov and Marcin Pilipczuk and Michal Pilipczuk and Pawel Rzazewski},
      year={2021},
      eprint={2007.11402},
      archivePrefix={arXiv},
      primaryClass={cs.DS},
      url={https://arxiv.org/abs/2007.11402}, 
}

@misc{paesani2022feedbackvertexsetcycle,
      title={Feedback Vertex Set and Even Cycle Transversal for H-Free Graphs: Finding Large Block Graphs}, 
      author={Giacomo Paesani and Daniël Paulusma and Paweł Rzążewski},
      year={2022},
      eprint={2105.02736},
      archivePrefix={arXiv},
      primaryClass={cs.DS},
      url={https://arxiv.org/abs/2105.02736}, 
}

@misc{bonnet2026qptasmwisfindinglarge,
      title={QPTAS for MWIS and finding large sparse induced subgraphs in graphs with few independent long holes}, 
      author={Édouard Bonnet and Jadwiga Czyżewska and Tomáš Masařík and Marcin Pilipczuk and Paweł Rzążewski},
      year={2026},
      eprint={2602.18317},
      archivePrefix={arXiv},
      primaryClass={cs.DS},
      url={https://arxiv.org/abs/2602.18317}, 
}

\end{document}